\documentclass[12pt]{amsart}
\usepackage{amscd,amssymb}
\usepackage[arrow,matrix]{xy}
\usepackage[plainpages,backref,urlcolor=blue]{hyperref}

\theoremstyle{plain}
\newtheorem{thm}[subsection]{Theorem}
\newtheorem{lem}[subsection]{Lemma}
\newtheorem{prop}[subsection]{Proposition}
\newtheorem{cor}[subsection]{Corollary}
\newtheorem{conj}[subsection]{Conjecture}

\theoremstyle{definition}
\newtheorem{rk}[subsection]{Remark}
\newtheorem{definition}[subsection]{Definition}
\newtheorem{ex}[subsection]{Example}
\newtheorem{question}[subsection]{Question}

\numberwithin{equation}{section}
\newcommand{\F}{{\mathcal F}}

\newcommand{\CC}{{\mathcal C}}

\newcommand{\NN}{{\mathcal N}}

\newcommand{\Q}{\mathbb{Q}}

\newcommand{\C}{\mathbb{C}}

\newcommand{\PP}{\mathbb{P}}

\DeclareMathOperator{\coker}{coker}

\DeclareMathOperator{\defect}{def}
\DeclareMathOperator{\codim}{codim}

\begin{document}

\title [On the syzygies and Alexander polynomials of nodal hypersurfaces]
{On the syzygies and Alexander polynomials of nodal hypersurfaces }

\author[Alexandru Dimca]{Alexandru Dimca$^1$}
\address{Institut Universitaire de France et Laboratoire J.A. Dieudonn\'e, UMR du CNRS 6621,
                 Universit\'e de Nice Sophia-Antipolis,
                 Parc Valrose,
                 06108 Nice Cedex 02,
                 France}
\email{dimca@unice.fr}

\author[Gabriel Sticlaru]{Gabriel Sticlaru}
\address{Faculty of Mathematics and Informatics,
Ovidius University,
Bd. Mamaia 124, 900527 Constanta,
Romania}
\email{gabrielsticlaru@yahoo.com }
\thanks{$^1$ Partially supported by the  ANR-08-BLAN-0317-02 (SEDIGA)}

\subjclass[2000]{Primary  14C30, 13D40; Secondary  32S35, 13D02}

\keywords{nodal hypersurfaces, Milnor algebra, syzygies, Alexander polynomials, Chebyshev polynomials, mixed Hodge structure, pole order filtration}

\begin{abstract} We give sharp lower bounds for the degree of the syzygies involving the partial derivatives of a homogeneous polynomial defining a nodal hypersurface. The result gives information on the position of the singularities of a nodal hypersurface expressed in terms of defects or superabundances.

The case of Chebyshev hypersurfaces is considered as a test for this result and leads to a potentially infinite family of nodal hypersurfaces having  nontrivial Alexander polynomials.

\end{abstract}

\maketitle


\section{Introduction} \label{sec:intro}

Let $S=\C[x_0,...,x_n]$ be the graded ring of polynomials in $x_0,,...,x_n$ with complex coefficients and denote by $S_r$ the vector space of homogeneous polynomials in $S$ of degree $r$. 
For any polynomial $f \in S_r$, we define the {\it Jacobian ideal} $J_f \subset S$ as the ideal spanned by the partial derivatives $f_0,...,f_n$ of $f$ with respect to $x_0,...,x_n$ 
and the corresponding graded {\it Milnor} (or {\it Jacobian}) {\it algebra} by
\begin{equation} 
\label{eq1}
M(f)=S/J_f.
\end{equation}
The study of such Milnor algebras is related to the singularities of the corresponding projective hypersurface $D:f=0$, see \cite{CD}, as well as to the mixed Hodge theory of 
the hypersurface $D$ and of its complement $U=\PP^n \setminus D$, see the foundational article by Griffiths \cite{Gr} and also \cite{Dc}, \cite{DS2}, \cite{DSW}, \cite{DSt2}.

The Milnor algebra $M(f)$ can be seen (up to a twist in grading) as the top cohomology of the Koszul complex $K^*(f)$ of the partial derivatives $f_0,...,f_n$ in $S$, see  \cite{CD} or  \cite{D1}, Chapter 6. As such, it is related to certain natural $E_1$-spectral sequences associated to the pole order filtration and converging to the cohomology of the complement $U$ introduced in \cite{Dc}, discussed in detail in \cite{D1}, Chapter 6 and reconsidered recently in \cite{DSt2}.

In the second section we study one of these spectral sequences for nodal hypersurfaces, using a key result by M. Saito telling when the Hodge filtration coincide to the pole order filtration on the cohomology groups $H^*(U)$. This study gives sharp lower bounds for the degree of the syzygies involving the partial derivatives of a homogeneous polynomial defining a nodal hypersurface,
extending the result proved in the curve case in Theorem 4.1 in \cite{DSt2} to arbitrary dimension.
In the curve case, see also \cite{RK} and \cite{EU}.

In the third section we consider the special case of Chebyshev hypersurfaces, which are classical examples of nodal hypersurfaces with many singularities. They were introduced by S. V. Chmutov to construct complex projective hypersurfaces with a large number of nodes, i.e. $A_1$-singularities, see \cite{AGV}, volume 2, p. 419 and \cite{Chm}.

In the final section,  we show  that on one hand the lower bounds obtained in the general case are best possible for curves and 3-dimensional Chebyshev hypersurfaces of degree $\leq 20$ (and probably for all odd dimensional  Chebyshev hypersurfaces, see Conjecture \ref{conj0}) , and on the other hand we give some topological applications, by computing the Alexander polynomials of Chebyshev hypersurfaces of dimension $2$ and $3$ and degree $d \leq 20$.

 The Alexander polynomials of singular hypersurfaces were introduced by A. Libgober \cite{Li1}, \cite{Li2} and are very subtle invariants of the topology of the complement $U$. However the number of classes of hypersurfaces where these Alexander polynomials are not trivial is rather limited, and this explains the interest of our new examples.

To end this Introduction, we recall the following notions, introduced in \cite{DSt2}.
\begin{definition}
\label{def}

For a hypersurface $D:f=0$ with isolated singularities we introduce three integers, as follows:

\noindent (i) the {\it coincidence threshold} $ct(D)$ defined as
$$ct(D)=\max \{q~~:~~\dim M(f)_k=\dim M(f_s)_k \text{ for all } k \leq q\},$$
with $f_s$  a homogeneous polynomial in $S$ of degree $d=\deg f$ such that $D_s:f_s=0$ is a smooth hypersurface in $\PP^n$.

\noindent (ii) the {\it stability threshold} $st(D)$ defined as
$$st(D)=\min \{q~~:~~\dim M(f)_k=\tau(D) \text{ for all } k \geq q\}$$
where $\tau(D)$ is the total Tjurina number of $D$, i.e. the sum of all the Tjurina numbers of the singularities of $D$.

\noindent (iii) the {\it minimal degree of a nontrivial syzygy} $mdr(D)$ defined as
$$mdr(D)=\min \{q~~:~~ H^n(K^*(f))_{q+n}\ne 0\}$$
where $K^*(f)$ is the Koszul complex of $f_0,...,f_n$ with the natural grading defined in \cite{DSt2}.

\end{definition}

 Moreover
it is easy to see that one has 
\begin{equation} 
\label{REL}
ct(D)=mdr(D)+d-2.
\end{equation} 

Recall also that, for a finite set of points $\NN \subset \PP^n$,
we denote by 
$$\defect S_m(\NN)=|\NN| - \codim \{h \in S_m~~|~~ h(a)=0 \text{ for any } a \in \NN\},$$ 
the {\it defect (or superabundance) of the linear system of polynomials in $S_m$ vanishing at the points in $\NN$}, see \cite{D1}, p. 207. This positive integer is called the {\it failure of $\NN$ to impose independent conditions on homogeneous polynomials of degree $m$} in \cite{EGH}.

When $D$ is a degree $d$ nodal hypersurface in $\PP^n$, with $\NN$ as singular set, it follows from Theorem 1.5 in \cite{DSt2} that one has 
\begin{equation} 
\label{def2}
 \defect S_k(\NN)\ne 0 \text{ for } k <T-ct(D) \text{ and } \defect S_k(\NN)= 0 \text{ for } k \geq T-ct(D)
\end{equation}
and also 
\begin{equation} 
\label{def3}
 \defect S_k(\NN)=|\NN|-\dim S_k \text{ for } k \leq T-st(D) 
\end{equation}
where $T=(n+1)(d-2).$
\bigskip
 
Note that computing the Hilbert-Poincar\'e series of the Milnor algebra $M(f)$ using an appropriate software is much easier than computing the defects $\defect S_k(\NN)$, because the Jacobian ideal comes with a given set of $(n+1)$ generators $f_0,...,f_n$, while the ideal $I$ of polynomials vanishing on $\NN$ has not such a given generating set.
However, it is the defects $\defect S_k(\NN)$, who describe the position of the singularities of $D$ in $\PP^n$ and which occur in many geometric problems, see for instance Theorem \ref{thmC} below.

Numerical experiments with the CoCoA package \cite{Co} and the Singular package \cite{Sing} have played a key role in the completion of this work.

\section{The spectral sequence and the syzygies of  nodal hypersurfaces} \label{sec2}

Let $D:f=0$ be a nodal hypersurface in $\PP^n$ of degree $d$. 

We consider first the case when $n=2n_1+1 \geq 3$ is odd. 
Then $D$ is a $\Q$-homology manifold
satisfying $b_j(D)=b_j(D_s)$ for $j \ne n-1$,  and the middle Betti number $b_{n-1}(D)$ is computable, e.g. using
the formula $b_{n-1}(D)=b_{n-1}(D_s)-n(D)$, 
where $n(D)=\tau (D)$ is the cardinal of the set $\NN$ of nodes of $D$. It follows that
the complement $U$ has at most two non-zero cohomology groups.
The first of them, $H^0(U)$ is $1$-dimensional and of Hodge type $(0,0)$, so nothing interesting here.
The second one, $H^n(U)$, is dual to $H^n_c(U)(-n)$ and $H^n_c(U)$ is isomorphic to $\coker (H^{n-1}(\PP^n) \to H^{n-1}(D))$, the morphism being induced by the inclusion $i:D \to \PP^n$. 

It follows that the mixed Hodge structure (for short MHS) on $H^n(U)$ is pure of weight $n+1$
with  
$$h^{p,q}(H^n(U))=h^{p-1,q-1}(D_s),$$
for $p+q=n+1=2n_1+2$, $p\ne q$, and 
$$h^{n_1+1,n_1+1}(H^n(U))=h^{n_1,n_1}(D)-1=h^{n_1,n_1}(D_s)-n(S)-1.$$
In particular, we have $P^1H^n(U)=F^1H^n(U)=H^n(U)$.

\medskip

But we have much more than this. Let $\alpha_D=\frac{n}{2}$. Then Corollary (0.12) in M. Saito \cite{Sa1}, or even better, the formula (1.1.3) in \cite{DSW}, imply that
\begin{equation} 
\label{FP1}
F^sH^n(U)=P^sH^n(U) \text{ for } s \geq n-\alpha_D+1
\end{equation}
i.e. for $s \geq n_1+2$.

\medskip

Now we look at the nonzero terms in the $E_1$-term of the spectral sequence $E_r^{p,q}(f)$
introduced in \cite{DSt2}, Proposition 2.2. Since $D$ has only isolated singularities,
these terms are sitting on two lines, given by $L: p+q=n$ and $L': p+q=n-1$.

\medskip

We look first at the terms on the line $L$. The term 
$E_1^{n-q,q}(f)=H^{n+1}(K^*(f))_{(q+1)d}$ is isomorphic as a $\C$-vector space to $M(f)_{(q+1)d-n-1}$, see Proposition (2.2) in \cite{DSt2}.

The corresponding limit term $E^{n-q,q}_{\infty}(f)=Gr_P^{n-q}H^{n}(U)$ coincides to $Gr_F^{n-q}H^{n}(U)$ in view of \eqref{FP1}, for $q=0,...,n_1-1=\left[ \frac{n}{2} \right]-1$,
where $\left[y \right]$ denotes the integral part of the real number $y$.

\medskip

On the other hand, Theorem 2.2 in \cite{DSW} yields the following isomorphism of $\C$-vector spaces
\begin{equation} 
\label{FP2}
Gr_F^{n-q}H^{n}(U)= M(f)_{(q+1)d-n-1}  \text{ for } q < \left[ \frac{n}{2} \right].
\end{equation}
It follows that in this range we have in fact
\begin{equation} 
\label{FP3}
E_1^{n-q,q}(f)= E^{n-q,q}_{\infty}(f).
\end{equation}
Therefore all the differentials in the $E_1$-spectral sequence $E_r(f)$ arriving at terms $E_r^{n-q,q}(f)$ having $q < \left[ \frac{n}{2} \right]$ are trivial.

\medskip

We look now at the terms on the line $L'$. The term 
$E_1^{n-1-q,q}(f)$ is given by $H^{n}(K^*(f))_{(q+1)d}$ and the limit term $E_{\infty}^{n-1-q,q}(f)$
has to be zero, as $H^{n-1}(U)=0$ for $n$ odd. It follows that
\begin{equation} 
\label{SZ1}
H^{n}(K^*(f))_{qd}=0  \text{ for } q \leq \left[ \frac{n}{2} \right]=n_1.
\end{equation}

\bigskip

We discuss now the case $n=2n_1 \geq 2$ even. Then $D$ is no longer a $\Q$-homology manifold, but one knows that $b_j(D)=b_j(D_s)$ for $j \notin \{ n-1,n\}$, and the $n$-th Betti number $b_{n}(D)$ is computable in terms of defects of linear series, namely
\begin{equation} 
\label{betti1}
b_n(D)=\defect S_{n_1d-2n_1-1}(\NN)+1,
\end{equation}
see Theorem (6.4.5) on page 208 in \cite{D1}. Moreover, the proof implies that the group $H^n(D)$ is a pure Hodge structure of type $(n_1,n_1)$, and the same holds for $H^{n-1}(U)$.

We consider the hypersurface $\tilde D$ in $\PP^{n+1}$ given by the equation $\tilde f(x_0,...,x_{n+1})=f(x_0,...,x_n)+x_{n+1}^d=0$ and the complement $\tilde U= \PP^{n+1} \setminus \tilde D$. Now $\tilde D$ is a $\Q$-homology manifold, and as above we define
 $\alpha_{\tilde D}=\frac{n}{2}+\frac{1}{d}$. Then  the formula (1.1.3) in \cite{DSW}, imply that
\begin{equation} 
\label{FP21}
F^sH^{n+1}(\tilde U)=P^sH^{n+1}(\tilde U) \text{ for } s \geq n_1+2.
\end{equation}

\medskip

Now we look at the nonzero terms in the $E_1$-term of the spectral sequence $E_r^{p,q}(\tilde f)$.
They are sitting on two lines, given by $L: p+q=n+1$ and $L': p+q=n$.

\medskip

We look first at the terms on the line $L$. The term 
$E_1^{n+1-q,q}(\tilde f)=H^{n+2}(K^*(\tilde f))_{(q+1)d}$ is isomorphic as a $\C$-vector space to $M(\tilde f)_{(q+1)d-n-2}$, see Proposition (2.2) in \cite{DSt2}.

The  limit term $E^{n+1-q,q}_{\infty}(\tilde f)=Gr_P^{n+1-q}H^{n+1}(\tilde U)$ coincides to $Gr_F^{n+1-q}H^{n+1}(\tilde U)$ in view of \eqref{FP21}, for $q=0,...,n_1-1=\left[ \frac{n}{2} \right]-1.$

\medskip

On the other hand, Theorem 2.2 in \cite{DSW} yields the following isomorphism of $\C$-vector spaces
\begin{equation} 
\label{FP22}
Gr_F^{n+1-q}H^{n+1}(\tilde U)= M(\tilde f)_{(q+1)d-n-2}  \text{ for } q < \left[ \alpha_{\tilde D} \right]=n_1.
\end{equation}
It follows that in this range we have in fact
\begin{equation} 
\label{FP23}
E_1^{n+1-q,q}(\tilde f)= E^{n+1-q,q}_{\infty}(\tilde f).
\end{equation}
All the differentials in the $E_1$-spectral sequence $E_r(\tilde f)$ arriving at terms $E_r^{n+1-q,q}(\tilde f)$ having $q < \left[ \frac{n}{2} \right]$ should therefore be trivial.

\medskip

Let us look now at the terms on the line $L'$. The term 
$E_1^{n-q,q}(\tilde f)$ is given by $H^{n+1}(K^*(\tilde f))_{(q+1)d}$ and the limit term $E_{\infty}^{n-q,q}(\tilde f)$
has to be zero, as $H^{n}(\tilde U)=0$ for $n$ even. It follows that
\begin{equation} 
\label{SZ21}
H^{n+1}(K^*(\tilde f))_{qd}=0  \text{ for } q \leq \left[ \frac{n}{2} \right].
\end{equation}
Now a non-zero class $\omega \in H^n(K^*(f))_m$ gives rise to a non-zero class $\tilde \omega =\omega \wedge dx_{n+1} \in H^{n+1}(K^*(\tilde f))_{m+1}$. It follows that
\begin{equation} 
\label{SZ22}
H^{n}(K^*(f))_{qd-1}=0  \text{ for } q \leq \left[ \frac{n}{2} \right]=n_1.
\end{equation}

Now recall that if the coordinates $x_0,...,x_n$ are choosen such that the hyperplane at infinity $H_0:x_0=0$ is transversal to $D$, then the multiplication by $x_0$ induces an injection
$H^n(K^*(f))_{m-1} \to H^n(K^*(f))_{m}$ (the dual statement for the homology is part of Corollary 11
in \cite{CD}).

This yields  our first main result.

\begin{thm}
\label{thmA} Let $D:f=0$ be a nodal surface in $\PP^n$ of degree $d$.

\smallskip

\noindent (i) If $n=2n_1+1$ is odd, then $H^n(K^*(f))_{m}=0$ for any $m \leq n_1d$.

\smallskip

\noindent (ii) If $n=2n_1$ is even, then $H^n(K^*(f))_{m}=0$ for any $m \leq n_1d-1$.

\end{thm}

If a single formula is preferred, then  $H^n(K^*(f))_{m}=0$ for any 
$$m \leq \left[ \frac{n}{2} \right]d-\frac{1+(-1)^n}{2}.$$

Using the formula \eqref{REL} we get the following.

\begin{cor}
\label{corA} Let $D:f=0$ be a nodal surface in $\PP^n$ of degree $d$ and let $\NN$ denote its set of nodes.

\smallskip

\noindent (i) If $n=2n_1+1$ is odd, then $ct(D)\geq (n_1+1)d-n-1=\frac{T}{2}$ and $\defect S_k(\NN)=0$ for
$k \geq (n_1+1)d-n-1=\frac{T}{2}$.

\smallskip

\noindent (ii) If $n=2n_1$ is even, then $ct(D) \geq (n_1+1)d-n-2=\frac{T}{2}+(\frac{d-2}{2}) $ and $\defect S_k(\NN)=0$ for
$k \geq n_1d-n= \frac{T}{2}-(\frac{d-2}{2}) $.

\end{cor}

\begin{ex}
\label{ex:n=2} 

\noindent (i) In the plane curve case, we have $n=2$, hence $(ii)$ in the above Theorem becomes
$H^2(K^*(f))_{m}=0$ for any $m \leq d-1$, which is exactly the first claim in Theorem 4.1 in \cite{DSt2}. Moreover this bound is strict, since $H^2(K^*(f))_d\ne 0$ for reducible curves as shown in loc.cit..

\noindent (ii) Corollary \ref{corA} (i) is a generalization of Theorem 5.1 (iv) in \cite{DSt2}. However, the formula in Remark \ref{rkAP} suggests that one should have 
the sharper bound
 $$ct(D) \geq \frac{T}{2}+(\frac{d-2}{2}) $$
 for $n$ odd as well.

\end{ex}

\section{Chebyshev  hypersurfaces: stability and coincidence threshholds} \label{sec3}

The $d$-th Chebyshev polynomial $T_d(x)=\cos (d \arccos (x))$ has $d-1$ critical points, namely $\lambda_k=\cos (k\pi/d)$ for $k=1,...,d-1$. One has $T_d(\lambda_k)=(-1)^k$. 

It follows that, for $d=2m+1$ odd,  the critical values $\pm 1$ are both attained $m$ times.
When $d=2m$, the maximal critical value $ 1$ is attained $m-1$ times, while the minimal critical value $-1$ is attained $m$ times.

\medskip

Consider the hypersurface $C(n,d,k)$ in $\PP^n$ for $n \geq 2$ defined by the homogeneous equation $f(n,d,k)=0$,
where the polynomial $f(n,d,k)$ is the homogenization (using $x_0$) of the polynomial
\begin{equation} 
\label{eqC}
g(n,k,d)=T_d(x_1)+ \dots + T_d(x_n)+k.
\end{equation} 
It is easy to see that the hypersurface $C(n,d,k)$ is smooth, unless $k$ is an integer satisfying
$|k|\leq n$ and $n+k$ is even. If these two conditions are fulfilled, then the hypersurface $C(n,d,k)$ is nodal, and the number of nodes is
$$ \tau(C(n,d,k))={n \choose a}d_1^n$$
if $d=2d_1+1$, with $2a=n+k$, and
$$ \tau(C(n,d,k))={n \choose a}d_1^n(1-\frac{1}{d_1})^a$$
if $d=2d_1$, with $2a=n+k$.

\medskip

It follows that for $d$ odd the maximal number of nodes is obtained for $a=\left[ \frac{n}{2} \right]$.
When $n=2n_1$ is even, this implies that $k=0$, so in this case the {\it Chebyshev hypersurface} $\CC(n,d)$ corresponds to $k=0$. For $n$ odd, both values $k=\pm 1$ give the same number of nodes. We pick the value $k=1$, for the reason explained below.

For $d$ even, it is not clear for which $k$  the maximum of $ \tau(C(n,d,k))$ is attained.
However, one may show that for $d \geq n+2$, the maximum is again attained for  $a=\left[ \frac{n}{2} \right]$.
We will call in this case the {\it Chebyshev} hypersurface $\CC(n,d)$ the hypersurface corresponding to
$a=\left[ \frac{n}{2} \right]$, $k=0$ for $n$ even, and $k=1$ for $n$ odd. In the latter case, i.e. $d$ even and $n$ odd, the choice $k=-1$ gives a lower number of nodes and it is less interesting, see Remark \ref{rkAP} $(ii)$.

In conlusion, the (affine part of) Chebyshev hypersurface $\CC(n,d)$ is defined  by the affine equation
$$g(n,d)=T_d(x_1)+ \dots + T_d(x_n)=0$$
when $n$ is even, and by 
$$g(n,d)=T_d(x_1)+ \dots + T_d(x_n)+1=0$$
when $n$ is odd.

Let $\NN(n,d)$ be the set of nodes of the  Chebyshev hypersurfaces $\CC(n,d)$. We may consider this set as a subset of the affine space $\C^n \subset \PP^n$ (given by $x_0=1$). For $n=2n_1$ even, the set $\NN(n,d)$ is the set of points $a=(a_1,...,a_n)$ such that $n_1$ among the $a_j$'s are local minimum points for $T_d$, and the remaining $n_1$ coordinates are local maximum points for $T_d$.
When $n=2n_1+1$, we have a similar description, the number of coordinates equal to local minima being $n_1+1$, and those of local maxima being $n_1$.

Consider the evaluation map
$$ev(n,d)_{\leq r}: \C[x_1,...,x_n]_{\leq r} \to \F(\NN(n,d))$$
where $\C[x_1,...,x_n]_{\leq r}$ denotes the vector space of polynomials of degree at most $r$, $\F(\NN(n,d))$
denotes the vector space of $\C$-valued functions on the set $\NN(n,d)$, and a polynomial $h$ is mapped to the function sending $a \in \NN(n,d)$ to $h(a)\in \C$. Then we have the following partial generalization of Proposition 3.1 in \cite{DSt1}.

\begin{prop}
\label{propB}
The evaluation map $ev(n,d)_r$ is injective if and only if $r \leq d-3$.
\end{prop}

\proof
Note that Proposition 3.1 in \cite{DSt1} implies the claim for $n=2$
Suppose first that $r \leq d-3$ and that the claim holds for $n-1 \geq 2$.
To fix the ideas, assume that $n=2n_1+1$ is odd. Fix $a_1$ to be one of the local minimum points of $T_d$. The set of points in $\NN(n,d)$ having the first coordinate equal to $a_1$ can be identified to the set of nodes $\NN(n-1,d)$, sitting in the affine space with coordinates $x_2,....,x_n$.
If $h \in \ker ev(n,d)_r$, it follows that $h(a_1,-) \in \ker ev(n-1,d)_r$. By our induction hypothesis, this kernel is trivial, hence $h(a_1,-)=0$ in the polynomial ring $\C[x_2,...,x_n]$.
It follows that $h$ is divisible by the polynomial $x_1-a_1$. Hence, in this way we get  $N_1 \geq d/2$ linear factors of $h$ of the form $x_1-a_1$. Using all the coordinates, we'll get  $N_n \geq nd/2>d$
distinct linear factors for $h$. This implies that $h=0$.

In the case $n=2n_1$ even, we should take $a_1$ a local maximum point of $T_d$ and all the rest goes in the same way as above.

To complete the proof it is enough to produce a polynomial $h \in \ker ev(n,d)$ with $\deg h =d-2$.
For this, let $f(n,d)(x_0,x_1,...,x_n)$ be the polynomial obtained from $g(n,d)$ by homogenization,
in other words $f(n,d)=0$ is an equation for $\CC(n,d)$ in $\PP^n$.
We take
$$h(x_1,...,x_n)=\frac{\partial f(n,d)}{\partial x_0}(1,x_1,...,x_n).$$
This polynomial vanishes on $\NN(n,d)$ by definition, and has degree $d-2$ because in the Chebyshev polynomial $T_d(x)$ the monomial $x^{d-1}$ is missing.

\endproof

Consider now the homogeneous ideal $I \subset S$ corresponding to the polynomials vanishing on the node set $\NN(n,d)$. The above result is equivalent to $I_k=0$ for $k<d-2$ and $I_{d-2} \ne 0$.
It follows that the corresponding defect
$$\defect S_k(\NN(n,d))=|\NN(n,d)|-\dim S_k+\dim I_k$$
satisfies $\defect S_k(\NN(n,d))=|\NN(n,d)|-\dim S_k$ for $k<d-2$ and $\defect S_k(\NN(n,d))>|\NN(n,d)|-\dim S_k$ for $k=d-2$. Using Theorem 1.5 in \cite{DSt2} we get the following improvement of Corollary 9 in \cite{CD}.

\begin{cor}
\label{corB}
Let $\CC(n,d)$ be the Chebyshev hypersurface of degree $d$ in $\PP^n$. Then the corresponding stability threshhold
is given by 
$$ st(\CC(n,d))=T-(d-3)=n(d-2)+1.$$
Moreover, the number of nodes is given by 

\smallskip

\noindent (i) $\tau(\CC(n,d))={2n_1 \choose n_1}d_1^n$ if $n=2n_1$ is even and $d=2d_1+1$ is odd;

\smallskip

\noindent (ii) $\tau(\CC(n,d))={2n_1 \choose n_1}d_1^{n_1} (d_1-1)^{n_1} $ if $n=2n_1$ is even and $d=2d_1$ is even;

\smallskip

\noindent (iii) $\tau(\CC(n,d))={2n_1+1 \choose n_1}d_1^{n_1+1} (d_1-1)^{n_1} $ if $n=2n_1+1$ is odd and $d=2d_1$ is even;

\smallskip

\noindent (iii) $\tau(\CC(n,d))={2n_1+1 \choose n_1}d_1^n $ if $n=2n_1+1$ is odd and $d=2d_1+1$ is odd.

\end{cor}

Note that the Chebyshev hypersurfaces are therefore among the very few classes of singular hypersurfaces $D$ for which the exact value of the stability threshold $st(D)$ is known.

Concerning the coincidence thresholds for Chebyshev hypersurfaces we have the following.

\begin{conj}
\label{conj0}

If $n=2n_1$ is even, then $ct(\CC(n,d))=(n_1+1)d-n-2$, i.e. the bounds given in Theorem \ref{thmA} and Corollary \ref{corA} are best possible in this case.
\end{conj}

\begin{ex}
\label{ex:ct} 
For $n=1$, the above Conjecture holds, see Example \ref{ex:n=2}.
For $n=4$, if we compute with {\it Singular} the Hilbert-Poincar\'e series of the Milnor algebra $M(f(n,d))$ for small values of $d$, say $3 \leq d \leq 20$, we see that the coincidence threshold is given in this case by the expected formula $ct(\CC(4,d))=3d-6$. The same holds for $n=6$, when $ct(\CC(6,d))=4d-8$.

\end{ex}

\section{Alexander polynomials of nodal hypersurfaces} \label{sec4}

Let $D$ be a degree $d$ hypersurface in $\PP^n$, with $d \geq 2$ and $n\geq 1$, given by a reduced equation $ f(x)=0$. 
Consider the corresponding global Milnor fiber $F$ defined by $f(x)-1=0$ in $\C^{n+1}$ with monodromy action $h:F \to F$, $h(x)=\exp(2\pi i/d)\cdot x$.
When $D$ has only isolated singularities, it is known that $\tilde H^k(F,\C)=0$ for $k<n-1$.
The characteristic polynomial
\begin{equation} 
\label{AlexPoly}
\Delta_D(t)=\det(t\cdot I-h^*|H^{n-1}(F,\C))
\end{equation}
is called the Alexander polynomial of the hypersurface $D$, with the convention $\Delta_D(t)=1$ if $b_{n-1}(F)=0$. To get a nontrivial Alexander polynomial $\Delta_D(t) \ne 1$, the idea is to look at hypersurfaces having lots of singularities, but sometime this is not enough, see Remark \ref{rkAP} (ii). That is why the number of examples  with $\Delta_D(t) \ne 1$ is rather limited.

For nodal hypersurfaces, one has the following precise description, see Theorem (6.4.5) on page 208 in \cite{D1} and our formula \eqref{def2}.

\begin{thm}
\label{thmC}
Let $D$ be a nodal hypersurface in $\PP^n$ of degree $d$ and let $\NN$ its set of nodes. Then:

\smallskip

\noindent (i) $\Delta_D(t)=1$ if $nd$ is odd;

\smallskip

\noindent (ii) $\Delta_D(t)=[t+(-1)^{n+1}]^{\defect S_m(\NN)}$ if $nd$ is even, where $m=nd/2-n-1$.
In this second case, $\Delta_D(t) \ne 1$ if and only if 
$$ct(D)<\frac{nd}{2}+d-n-1=\frac{T}{2}+\frac{d}{2}.$$
\end{thm}

\begin{rk}
\label{rkAPN}

(i) For $n$ even, the lower bound given by Corollary \ref{corA} can be written as
$$\frac{nd}{2}+d-n-2 \leq ct(D).$$
It follows that in the case $n$ even, $\Delta_D(t) \ne 1$ if and only if 
$$ct(D)=\frac{nd}{2}+d-n-2.$$
This seems to be the case for all odd dimensional Chebyshev hypersurfaces, see Conjecture \ref{conj1},
Example \ref{exC} and Remark \ref{rkAP}. This explains our special attention devoted to this class of nodal hypersurfaces.

(ii) One may express the topological content of Theorem \ref{thmC} using only Betti numbers as follows.
If $n=2n_1$ is even, then one has
$$b_n(D)=1+\defect S_m(\NN)$$
where $m=nd/2-n-1$. It follows that $\dim H^n(K^*(f))_{n_1d}=b_n(D)-1$, i.e. $mrd(D) \geq n_1d-n$ and the number of nontrivial syzygies of the minimal expected degree $m=n_1d-n$ is determined topologically, exactly as in the case $n=2$ covered by Theorem 4.1 in \cite{DSt2}.
However, for $n>2$, we do not have explicit formulas for these syzygies.

For $n$ odd and $d$ even, let $D^2$ be the double cover of $\PP^n$ ramified along $D$. Then
$$b_{n+1}(D^2)=1+\defect S_m(\NN).$$
\end{rk}

\begin{ex}
\label{kummer}
Consider the {\it Kummer surface} $S$ in $\PP^3$ given by the affine equation
$$x^4 + y^4 + z^4 - y^2z^2 - z^2x^2 - x^2y^2 - x^2 - y^2 - z^2 + 1=0,$$
see for instance \cite{Hu}, p. 93.
This surface has the maximum number of nodes for a surface in $\PP^3$ of degree $4$, namely $16$ nodes. A direct computation with Singular yields
$$HP(M(f))(t)=1+4t+10t^2+16t^3+19t^4+16(t^5+t^6+…)$$
and hence for the Betti number of the associated $3$-fold $D^2$ we get
$$b_{4}(D^2)=1+\defect S_2(\NN)=7.$$
\end{ex}

Using Theorem \ref{thmC}
we may get a (potentially infinite) family of Chebyshev hypersurfaces in $\PP^n$ for $n=3$ and $n=4$
with rather large Alexander polynomials. We offer the following.

\begin{conj}
\label{conj1}

\noindent (i) Let $\CC(3,d)$ be the Chebyshev surface of even degree $d=2d_1$ in $\PP^3$. Then
$$\defect S_{3d_1-4}(\NN(3,d))=3(d_1-1).$$
\noindent (ii) Let $\CC(4,d)$ be the Chebyshev $3$-folds of degree $d$ in $\PP^4$. Then
$$\defect S_{2d-5}(\NN(4,d))= \left[ \frac{d-1}{2} \right] ( 3\left[ \frac{d-1}{2} \right] -1)  .$$

\end{conj}

\begin{ex}
\label{exC}
Conjecture \ref{conj1} holds for all degrees $d$ satisfying $3 \leq d \leq 20$.
\end{ex}

 As explained in the Example \ref{ex:ct}, for $n=4$ the bounds given in Corollary \ref{corA} are best possible. In particular we have $\defect S_{k}(\NN(4,d))=0$ for $k>2d-5$,
but $\defect S_{2d-5}(\NN(4,d))>0$. The actual values are computed using the Singular software.

In the case $n=3$, the bounds given in Corollary \ref{corA} are not optimal.  But we compute directly the Hilbert-Poincar\'e series of the Milnor algebra $M(f(n,d))$ in this range and use the relations between the coefficients and the defects $\defect S_{k}(\NN(3,d))$ given by Theorem 1.5 in
\cite{DSt2}. 

\endproof

\begin{rk}
\label{rkAP}

\noindent (i) It is interesting to note that in both cases above, one has
$$ct(D)=\frac{nd}{2}+d-n-2,$$
i.e. the inequality in Theorem \ref{thmC} is very tight in the case of Chebyshev hypersurfaces.

\smallskip

\noindent (ii)
The feeling that the construction of nodal hypersurfaces $D$ with $\Delta_D(t) \ne 1$ is not easy is reflected by the fact that in the case $n=3$ all the even degree surfaces obtained from the Chebyshev surface by changing $k=1$ into $k=-1$ in the equation \eqref{eqC} have $\Delta_D(t)= 1$, in spite of being  nodal surfaces with lots of nodes.

\smallskip

\noindent (iii) Similar results hold for the Alexander polynomials of Chebyshev hypersurfaces of dimension $\geq 4$. We leave the interested reader to find the exact statements in each case. To do this, he has to compute the Hilbert-Poincar\'e series $HP(M(f))$ using a computer algebra software, check whether
$$a=\frac{nd}{2}+d-n-1-ct(D)>0,$$
and then, in view of Theorem 1.5 in \cite{DSt2}, compute the difference
$$\defect S_m(\NN)=\dim M(f)_{ct(D)+a}-\dim M(f_s)_{ct(D)+a}.$$
Here $m=nd/2-n-1$ and $D_s:f_s=0$ is any smooth hypersurface in $\PP^n$ of degree $d$.

\end{rk}

\end{document}